\documentclass[11pt,a4paper]{article} 
\usepackage[cp866]{inputenc}
\usepackage[english,russian]{babel}
\usepackage{amsmath}
\usepackage{amsfonts}
\usepackage{longtable}
\usepackage{amssymb}
\usepackage{graphics}
\vspace{2.5cm}

\begin{document}
УДК 517.172.4 \vspace{0.5cm}

\begin{center} {\bf Прудников И.М.} \end{center} \vspace{0.5cm}
\begin{center}
{\large \bf К ВОПРОСУ О ПРЕДСТАВИМОСТИ ФУНКЦИИ ДВУХ ПЕРЕМЕННЫХ В
ВИДЕ РАЗНОСТИ ВЫПУКЛЫХ ФУНКЦИЙ (исправленный вариант)}
\end{center}
\vspace{0.5cm}

В статье приведены исправленные необходимые и достаточные условия
представимости произвольной липшицевой функции двух переменных в
виде разности выпуклых функций \cite{proudconvex2}. Дана также
геометрическая интерпретация этих условий. Приведен алгоритм
такого представления, результатом которого есть последовательность
равномерно сходящихся выпуклых функций.

\noindent

{\bf Ключевые слова:} липшицевые функции, выпуклые функции,
вариация функции, кривизна кривой, поворот кривой.

\vspace{1cm}
\section{Введение}

Эта проблема была впервые сформулирована академиком
А.Д.Александровым в статье \cite{aleksandrov1} и исследована
многими российскими и зарубежными  математиками
\cite{aleksandrov2} - \cite{Hartman}. Решение этой проблемы
интересно как для геометров, так и для математиков, занимающихся
оптимизацией, например, для построения квазидифференциального
исчисления \cite{demvas}.

Необходимые и достаточные условия представимости функции одной
переменной в виде разности выпуклых, т.е. условия. когда функция
является ПРВ функцией, хорошо известны. Эти условия могут быть
записаны в следующем виде.

Пусть $x \rightarrow f(x): [a,b] \rightarrow \mathbb{R}$ -
произвольная липшицевая функция. Известно, что множество $N_f$,
где функция $ f(\cdot)$ дифференцируемая, есть множество полной
меры на [a,b]. Для того, чтобы функция $f(\cdot)$ была представима
в виде разности выпуклых функций, необходимо и достаточно, чтобы
выполнялось условие
$$
         \vee (f'; a,b) < \infty ,
$$
где производные вычисляются там, где они существуют. Символ $\vee$
означает вариацию функции $ f'$ на отрезке [a,b].

В той же статье \cite{aleksandrov1}  А.Д.Александров задает вопрос
о представимости функции в виде разности выпуклых, если она
является таковой для любой прямой в области определения. Ответ на
этот вопрос отрицательный (см. \cite{hiriarturruty}, \cite{vesely}
).

Согласно терминологии А.Д.Александрова под многогранной
кусочно-линейной функцией с конечным числом граней будем понимать
такую функцию, график которой состоит из конечного числа частей
плоскостей, которые называются гранями.

В статье  \cite{proudconvex1} даны необходимые и достаточные
условия представимости произвольной липшицевой положительно
однородной (п.о.) функции трех переменных в виде разности выпуклых
функций. Результат может быть распространен на положительно
однородные функции $m$-ой степени. Теперь откажемся от условия
положительной однородности и будем рассматривать произвольную
липшицевую функцию $f(\cdot)$ с константой Липшица $L$ от двух
переменных $(x,y) \rightarrow f(x,y):D \rightarrow \mathbb{R},$
где $D$ есть выпуклое открытое ограниченное множество в
$\mathbb{R}^2$, так что его замыкание $\bar{D}$ - компакт.
Приведем алгоритм такого представления и найдем необходимые и
достаточные условия сходимости построенной последовательности
функций.

Дадим определение звездной области.

{\bf Определение1}. {\em Область ${D}$ называется звездной, если
существует точка $O \in \int {D} $, для которой отрезок,
соединяющий точку $O$ с любой точкой области ${D}$, целиком
принадлежит области $D$.}

Пусть ${\wp}(D)$ - класс кривых на плоскости  $XOY$ в множестве
$D$, ограничивающих звездные области, замыкание которых компактно.
Параметризуем кривые $r \in {\wp}(D)$ естественным образом, т.е.
параметр $\tau$ точки $M$ на кривой $r(\cdot)$ равен длине кривой
между $M$ и начальной точкой. Обозначим такую кривую как $r(t), t
\in [0,T_r]$.

Заметим, что множеству кривых ${\wp}(D)$ принадлежит класс кривых
на плоскости  $XOY$ в множестве $D$, параметризованных
естественным образом и ограничивающих выпуклые компактные
множества на множестве $D$. Обозначим этот класс кривых через
$\varrho(D)$.

С помощью кривых  $r \in \wp (D)$ необходимые и достаточные
условия представимости функции $f(\cdot)$ в виде разности выпуклых
функций  могут быть записаны в следующем виде.

{\bf Теорема 1.} {\em Для того, чтобы липшицевая функция $z
\rightarrow f(z):D \rightarrow \mathbb{R}$ была представима в виде
разности выпуклых функций (была ПРВ функцией), необходимо и
достаточно, чтобы для любой кривой $r \in \wp (D)$ и любой ее
системы подмножеств выполнялось неравенство
$$
 (\exists c_1(D,f), c_2(D,f)>0) (\forall r \in \wp(D)) ,
$$
\begin{equation}
 \vee (\Phi'; 0,T_r)<  c_1(D,f) + c_2(D,f) \vee (r';
 0,T_r), \label{difconv3}
\end{equation} где $\Phi(t)=f(r(t)) \;\;\;\; \forall t \in
[0,T_r]$. \label{difconvexthm1} }

{\bf Замечание 1}. {\em  Условие теоремы означает следующее.

Существуют константы $ c_1(D,f), c_2(D,f)>0 $, что для любой $r
\in \wp(D)$ и любой системы подмножеств $[T_i, T_{i+1} ] \subset
[0, T_r] $, $i \in 1:m,$
$$
 \sum_1^m \vee (\Phi'; T_i,T_{i+1})<  c_1(D,f) + c_2(D,f) \sum_1^m\vee (r';
 T_i,T_{i+1}).
$$
Необходимость рассмотрения системы подмножеств для кривой $r \in
\wp (D)$, для которой выполняется неравенство (\ref{difconv3}),
объясняется тем, что вариация $\vee (r'; 0,T_r)$ может быть
неограниченной для всей кривой $r(\cdot)$.}

{\bf Замечание 2}. {\em  Впервые подобные условия для вариации
функции на кривых были приведены  в \cite{proudconvex1}. }

Доказательство основано на специальном алгоритме представления
функции $f(\cdot)$ в виде разности выпуклых функций. В результате
получаем конечные или бесконечные последовательности выпуклых
функций, равномерно сходящиеся на $D$  к выпуклым функциям,
разность которых есть исходная функция $f(\cdot)$, если условия
теоремы \ref{difconvexthm1} выполняются.

Ниже приведен алгоритм представления функции $f(\cdot)$ в виде
разности выпуклых функций и доказана его сходимость, если условия
теоремы \ref{difconvexthm1} выполняются.

Для представления функции $f(\cdot)$ в виде разности выпуклых
функций будем использовать две операции, в результате которых
получаем конечное или счетное число выпуклых многогранных
кусочно-линейных функций, определенных на $D$.

{\em Первая операция} - это приближение функции $f(\cdot)$
многогранной кусочно-линейной функцией $f_n(\cdot)$ с конечным
числом граней.

{\em Вторая операция} - это представление функции $f_n(\cdot)$ в
виде разности выпуклых многогранных кусочно-линейных функций
$f_{1,n}(\cdot):D \rightarrow \mathbb{R}$ и $f_{2,n}(\cdot):D
\rightarrow \mathbb{R}$ согласно алгоритму, описанному ниже.

Далее доказывается, что, если выполняются условия теоремы
\ref{difconvexthm1}, то из последовательностей $f_{1,n}(\cdot)-
c_{1,n}$ и $f_{2,n}(\cdot) - c_{1,n}$, где $c_{1,n}=f_{1,n}(a)$,
$a - $ произвольная внутренняя точка области $D$, можно выбрать
сходящиеся подпоследовательности.

Когда условия теоремы выполняются, то, как будет показано,
вариация производной вдоль любого отрезка множества $D$ выпуклых
функций $f_{1,n}(\cdot)$ и $f_{2,n}(\cdot)$ ограничена сверху
константой, зависящей от $D$ и $f$.

Метод представления конечной многогранной функции в виде разности
выпуклых подобен методу, использованному А.Д.Александровым в
\cite{aleksandrov1} при исследовании возможности представления
специального вида функций в виде разности выпуклых.

\vspace{1cm}

\section {Доказательство теоремы}

Начнем доказательство теоремы с описания алгоритма.

\vspace{0.5cm}

{\bf ОПИСАНИЕ АЛГОРИТМА}

\vspace{0,5cm}

1. Производим достаточно равномерную триангуляцию области $D$ и
строим по каждому треугольнику линейную функцию, значения которой
равны значениям функции $f(\cdot)$ в вершинах треугольника.
Функцию с получившимся графиком обозначим через $f_n(\cdot):D
\rightarrow \mathbb{R}$, где $n$ равно числу треугольников, на
которые мы разбиваем область $D$.

2. Представляем функцию $f_n(\cdot)$ в виде разности выпуклых
согласно алгоритму, как это описано ниже.

Предварительно введем понятие двугранного угла. Будем понимать под
двугранным углом функцию, определенную на $D$, график которой
состоит из полуплоскостей с общей граничной прямой, называемого
ребром двугранного угла.

Рассмотрим все выпуклые двугранные углы, части графиков которых
принадлежат графику функции $f_n(\cdot)$. Определяем эти
двугранные углы на всей области $D$. Просуммируем все такие
выпуклые двугранные углы. В итоге получим выпуклую многогранную
функцию $f_{1,n}(\cdot):D \rightarrow \mathbb{R}$. Доказывается
\cite{aleksandrov1}, что разность \begin{equation}
 f_{1,n}(\cdot) - f_n(\cdot) = f_{2,n}(\cdot)
\label{difconv1} \end{equation} есть также выпуклая многогранная
функция.

Действительно, для доказательства достаточно показать, что все
двугранные углы, части графиков которых принадлежат графику
функции $f_{1,n}(\cdot) - f_{n}(\cdot), - $ выпуклые. Для этого
покажем, что любая точка, лежащая на проекции ребра произвольного
двугранного угла функции $f_{1,n}(\cdot) - f_{n}(\cdot),$ имеет
малую окрестность, где функция $f_{1,n}(\cdot) - f_{n}(\cdot)$
выпукла.

Если берем точку, в малой окрестности которой функция
$f_{n}(\cdot)$ линейна, то локальная выпуклость разности
$f_{1,n}(\cdot) - f_{n}(\cdot)$ очевидна. Пусть берем точку,
лежащую на проекции на плоскость ребра выпуклого двугранного угла,
часть графика которого принадлежит графику функции $f_{n}(\cdot)$.
Поскольку согласно алгоритму этот же двугранный угол входит в
сумму выпуклых двугранных углов, образующих функцию
$f_{1,n}(\cdot)$, то опять разность $f_{1,n}(\cdot) -
f_{n}(\cdot)$ будет локально выпуклой в окрестности
рассматриваемой точки. Если же точка лежит на проекции ребра
вогнутого двугранного угла, часть графика которого принадлежит
графику функции $f_{n}(\cdot)$, то $-f_{n}(\cdot)-$ локально
выпукла в окрестности этой точки, а поэтому разность
$f_{1,n}(\cdot) - f_{n}(\cdot)$ снова локально выпукла в той же
окрестности. Из локальной выпуклости всех двугранных углов функции
$f_{1,n}(\cdot) - f_{n}(\cdot)$ следует ее выпуклость на всем
множестве $D$.

Покажем, что при выполнении теоремы \ref{difconvexthm1} из
последовательности функций $f_{1,n}(\cdot) - c_{1,n} $ можно
выделить подпоследовательность, равномерно сходящуюся  на $D$ к
выпуклой функции $f_1(\cdot)$ при $n \rightarrow +\infty$. Тогда
из (\ref{difconv1}) будет следует, что подпоследовательность
функций $f_{2,n}(\cdot)-c_{1,n}$ также равномерно сходится к
выпуклой функции $f_2(\cdot)$. Для функций $f_1(\cdot)$ и
$f_2(\cdot)$ верно равенство \begin{equation}
 f_{1}(\cdot) - f_2(\cdot) = f(\cdot).
\label{difconv2} \end{equation}

Начнем доказательство с одномерного случая, когда $D=[a,b] \subset
\mathbb{R}$.

Приблизим функцию $f(\cdot)$ кусочно-линейной функцией
$f_n(\cdot)$  с любой степенью точности. На первом шаге выделяем
все выпуклые двугранные углы, части графиков которых принадлежат
графику функции $f_n(\cdot)$. Распространяем их на весь отрезок
$[a,b]$ и просуммируем. В итоге получим выпуклую кусочно-линейную
функцию $f_{1,n}(\cdot):[a,b] \rightarrow \mathbb{R}$. Согласно
сказанному выше разность $f_{1,n}(\cdot) - f_{n}(\cdot)$ есть
снова выпуклая кусочно-линейная функция на $[a,b]$.

Покажем, что вариация производных функций $f_{1,n}(\cdot)$ и
$f_{2,n}(\cdot)$ на отрезке $[a,b]$ ограничена сверху той же
константой $c$, что вариация производной функции $f(\cdot)$, т..е.
$$
\vee (f_{1,n}' ; a,b) \leq  c.
$$
Последнее следует из цепочки неравенств
$$
\vee (f_{1,n}' ; a,b) \leq \vee (f_n' ; a,b) \leq \vee (f' ; a,b)
\leq c.
$$
Но тогда из $f_{1,n}(\cdot)$ можно вычесть константу
$c_{1,n}=f_{1,n}(a)$, $ a \in \mbox{int} \, D$, чтобы функции
$f_{1,n}(\cdot)-c_{1,n}$ были ограниченными на отрезке $[a,b]$ в
совокупности по $n$, т.е. равностепенно ограниченными. Из оценки
для вариации производной, не зависящей от $n$, следует
равностепенная непрерывность функций $f_{1,n}(\cdot)-c_{1,n}$. Из
теоремы Арцела  получим, что из последовательности выпуклых
функций $f_{1,n}(\cdot) - c_{1,n}$ можно выделить
подпоследовательность функций $f_{1,n_k}(\cdot) - c_{1,n_k}$,
которая сходится равномерно при $n_k \rightarrow \infty$ к
некоторой выпуклой на $[a,b]$ функции $f_{1}(\cdot)$.
Соответственно, последовательность функций $f_{2,n_k}(\cdot) -
c_{1,n_k}$ также равномерно на $[a,b]$ сходится при $n_k
\rightarrow \infty$ к некоторой выпуклой на $[a,b]$ функции
$f_{2}(\cdot)$. В итоге будем иметь
$$
f(\cdot)=f_1(\cdot)- f_2(\cdot).
$$

Перейдем к двумерному случаю и покажем, что тот же алгоритм
приводит к паре выпуклых функций на $D$, разность которых есть
исходная функция $f(\cdot)$.

Возьмем произвольную кривую $r(\cdot) \in \wp(D).$Пусть
$$
    \Phi(t)=f(r(t)) \;\;\; \forall t \in [0,T_r].
$$
Покажем, что $\Phi(\cdot)$ - липшицевая с константой $L$.
Действительно, для любых $t_1,t_2 \in [0,T_r]$ имеем
$$
 \mid \Phi (t_1) - \Phi (t_2) \mid = \mid f(r(t_1)) - f(r(t_2)) \mid \le
 L \Vert  r(t_1) - r(t_2) \Vert  \le  L \mid t_1 - t_2 \mid.
$$
Поэтому \cite{kolmogorovfomin} $\Phi (\cdot)$ почти всюду (п.в.)
дифференцируемая на $[0,T_r].$ Множество точек дифференцируемости
функции $\Phi (\cdot)$ на $[0,T_r]$ обозначим через $N_r.$

Докажем, что если существуют константы $c_1(D, f), c_2(D,f) > 0$
такие, что для произвольной кривой $r(\cdot) \in \wp(D)$
выполняется условие теоремы 1, то из последовательностей функций
$f_{1,n}(\cdot)-c_{1,n}$, $f_{2,n}(\cdot)-c_{1,n}$,  можно выбрать
подпоследовательности, равномерно на $D$ сходящиеся к выпуклым
функциям $f_{1}(\cdot)$, $f_{2}(\cdot)$ соответственно, для
которых верно равенство (\ref{difconv2}).

Доказательство будем основывать на леммах, приведенных ниже.

{\bf Лемма 1.} {\em Для любой выпуклой п.о. степени 1 функции $q
\rightarrow \psi (q): \mathbb{R}^2 \rightarrow \mathbb{R}$ и
кривой $ r(\cdot) \in \wp(D)$, а также любого ее подмножества
верно неравенство
$$
    \vee (\Theta '; 0,T_r)  <  c_1(D, \psi) + c_2(D, \psi)
    \vee(r'; 0,T_r ),
$$
где $\Theta(t) = \psi(r(t))$  для всех $t \in [0,T_r], c_1(D,
\psi), c_2(D, \psi)$ - некоторые положительные константы.
\label{1lemdifconv} }

{\bf Доказательство.}  Без ограничения общности будем считать, что
$\psi(\cdot)$ есть гладкая функция на $\mathbb{R}^2 \backslash
\{0\}.$ Пусть
$$
 \psi (r(t)) = \mbox{max} \,_{v \in \partial \psi(0)} (v, r(t))=(v(t),r(t)), \;\;
v(t) \in \partial \psi(0),
$$
где $\partial \psi (0)$ - субдифференциал функции $\psi(\cdot)$ в
нуле. Будем также считать,  что $r(\cdot)$ - дифференцируемая
кривая по $t \in [0,T_r]$ .

Очевидно, что
$$
\psi'(r(t))=(v'(t),r(t))+(v(t),r'(t)).
$$
Так как $r(t)$ есть нормаль к границе множества $\partial \psi(0)$
в точке $v(t)$, то векторы $v'(t)$ и $r(t)$ перпендикулярны друг к
другу, а следовательно, $(v'(t),r(t))=0.$  Поскольку кривая
$r(\cdot)$ параметризована естественным образом, то $\Vert r'(t)
\Vert =1$ для любых $t \in [0,T_r]$.

Нетрудно проверить следующую цепочку неравенств
$$
\mid \psi' (r(t_1)) - \psi' (r(t_2)) \mid= \mid (v(t_1),r'(t_1)) -
(v(t_2),r'(t_2)) \mid = \mid (v(t_1)- v(t_2),r'(t_1)) +
$$ $$
(v(t_2),r'(t_1))- (v(t_2),r'(t_2)) \mid \leq  \Vert v(t_1) -
v(t_2) \Vert \; \Vert r'(t_1) \Vert + \Vert r'(t_1) - r'(t_2)
\Vert \; \Vert v(t_2) \Vert
$$
Отсюда следует, что
$$
   \vee (\Theta '; 0,T_r) < 2P(\partial \psi (0)) + L(D) \vee(r'; 0; T_r ),
$$
где $P(\psi (0))$- длина кривой, ограничивающей выпуклое
компактное множество $\partial \psi(0) \subset \mathbb{R}^2$ и
$L(D)$ - константа Липшица функции $\psi(\cdot).$

Пусть теперь $\psi (\cdot)$ - произвольная выпуклая ПО функция. С
любой степенью точности ее можно приблизить  на единичном круге
выпуклой ПО дифференцируемой на $\mathbb{R}^2 \backslash \{0\}$
функцией $\hat{\psi}(\cdot)$ так, чтобы в метрике Хаусдорфа
субдифференциалы в нуле этих функций отличались друг от друга как
угодно мало. Но тогда и длины кривых, ограничивающих их
субдифференциалы, будут отличаться друг от друга как угодно мало.
Также кривую $r(\cdot)$ можно приблизить дифференцируемой кривой
таким образом, чтобы их производные по $t$ в точках
дифференцируемости кривой $r(\cdot)$ отличались друг от друга по
норме на произвольно малую величину. Таким образом, любые конечные
суммы, используемые при вычислении вариаций  функций
$\Theta'(\cdot)$ и $\hat{\Theta}'(\cdot)$  для негладкого и
гладкого случая, могут быть сделаны за счет приближения как угодно
близкими друг к другу. Но поскольку вариацию функции
$\hat{\Theta}'(\cdot)$ можно ограничить сверху величиной,
зависящей только от множества $D$ и некоторых констант, то лемма 1
доказана. $\Box$

На основе этой леммы докажем утверждение.

{\bf Лемма 2.} {\em Пусть $(x,y) \rightarrow f_1(x,y):
\mathbb{R}^2 \rightarrow \mathbb{R}-$ непрерывная выпуклая функция
и $r(\cdot) \in \wp(D), t \in [0,T_r].$ Тогда существуют константы
$c_1(D, f_1), c_2(D,f_1)>0,$ что для кривой $r(\cdot)$ и любого ее
подмножества
\begin{equation}
   \vee (\Phi_1'; 0,T_r)  <  c_1(D, f_1) + c_2(D, f_1) \vee(r'; 0,T_r ),
\label{difconv4} \end{equation} где $\Phi_1(t) = f_1(r(t)), t \in
[0,T_r].$ \label{2lemdifconv} }

{\bf  Доказательство.} На начальном этапе будем считать, что
$f_1(\cdot,\cdot)$ дважды непрерывно дифференцируемая функция на
$D$, которая принимает неотрицательные значения и начало координат
$-$ ее точка минимума, а также, что $0=(0,0)$ принадлежит
внутренности выпуклой области на $\mathbb{R}^2$ с границей
$r(\cdot).$

Построим для функции $f_1(\cdot,\cdot)$ п.о. степени 1 функцию
$\psi(\cdot),$  которая на $r(\cdot)$ принимает значения, равные
$f_1(r(\cdot)).$ Покажем, что $\psi(\cdot)$ - выпуклая.

Рассмотрим функцию
$$
f_{\varepsilon}(x,y)=f_1(x,y)+\varepsilon ( \mid \mid x \mid
\mid^2 + \mid \mid y \mid \mid^2  ), \,\,\, \varepsilon>0.
$$
Разобьем отрезок $[ 0,T_r] $ точками  $\{t_i\} , i \in 1:J ,$ на
равные отрезки. Построим плоскости $\pi_i$ в $\mathbb{R}^3,$
проходящие соответственно через точки $(0,0,0),
(r(t_i),f_{\varepsilon}(r(t_i))),
(r(t_{i+1}),f_{\varepsilon}(r(t_{i+1})), i \in 1:J $. Части
плоскостей $\pi_i ,i \in 1:J$ , определенных в секторах,
образуемых векторами $(0,0), r(t_i), r(t_{i+1})$, определяют
график п.о. степени 1 многогранной функцию $(\psi_{\varepsilon})_J
(r(\cdot)).$ Будем понимать под двугранным углом функцию, график
которой состоит из полуплоскостей с общей граничной прямой,
включающих плоскости $\pi_i,$ построенные в соседних секторах.
Покажем, что все двугранные углы функции $(\psi_{\varepsilon})_J
(r(\cdot)),$ образуемые смежными плоскостями $\pi_i, i \in J,-$
выпуклые.

Поскольку всегда любую кривую $r(\cdot) \in \wp(D)$ можно
приблизить  с любой степенью точности гладкой кривой из $\wp(D),$
то без ограничения общности будем считать, что $r(\cdot)$ -
гладкая дифференцируемая кривая с производной $r'(\cdot).$

Под градиентом плоскости $\pi_i$ будем понимать градиент линейной
функции, график которой совпадает с плоскостью $\pi_i$. Обозначим
градиенты плоскостей $\pi_i$ и $\pi_{i+1}$ через $\nabla \pi_i$ и
$\nabla \pi_{i+1}$ соответственно.  Воспользуемся теоремой о
средней точке, согласно которой существует такая точка $t_m \in
[t_i, t_{i+1}],$ что
$$
           \partial f_{\varepsilon}(r(t_m))/ \partial e_i =
           (\nabla \pi_i, e_i),
$$
где
$$
e_i=(r(t_{i+1})-r(t_i))/ \mid \mid r(t_{i+1})-r(t_i) \mid \mid .
$$
Аналогично для плоскости $\pi_{i+1}$ и некоторой точки $t_c \in
[t_{i+1},  t_{i+2}]$ имеем
$$
           \partial f_{\varepsilon}(r(t_c))/ \partial e_{i+1} =
           (\nabla \pi_{i+1}, e_{i+1}),
$$
где
$$
e_{i+1}=(r(t_{i+2})-r(t_{i+1}))/ \mid \mid r(t_{i+2})-r(t_{i+1})
\mid \mid .
$$
Функция $f_{\varepsilon}(\cdot)$ сильно выпуклая, так как ее
матрица вторых частных производных положительно определенная.
Любая выпуклая функция имеет неубывающую производную по
направлению вдоль произвольного луча. Но для сильно выпуклой
функции производная по касательному направлению к кривой вида
$r(x_0,\tau, g)= x_0+\tau g +o_{\varepsilon}(\tau)$, $g \in
\mathbb{R}^n$,$ \tau >0$ в малой окрестности точки $x_0$ есть
возрастающая функция вдоль этой кривой. Поэтому для достаточно
большом $J$ и равномерном разбиении кривой $r(\cdot)$ точками
$t_i$ имеем
$$
\partial f_{\varepsilon}(r(t_m))/ \partial e_i <
\partial f_{\varepsilon}(r(t_c))/ \partial e_{i+1},
$$
или
$$
(\nabla \pi_i, e_i) < (\nabla \pi_{i+1}, e_{i+1}).
$$
Учтем также, что разность $\nabla \pi_{i+1} - \nabla \pi_i$
перпендикулярна вектору $r(t_{i+1}).$ Отсюда и из неравенства выше
следует, что двугранный угол $\pi_i, \pi_{i+1}$ - выпуклый. При $J
\rightarrow \infty$
$$ (\psi_{\varepsilon})_J(\cdot) \Rightarrow (\psi_{\varepsilon})(\cdot).$$
Так как точечный предел для выпуклых функций равносилен
равномерному пределу, то $\psi_{\varepsilon}(\cdot)$ - выпуклая
функция. Также $\psi_{\varepsilon}(\cdot) \Rightarrow \psi
(\cdot)$ при $\varepsilon \rightarrow +0,$ т.е. $\psi(\cdot)-$
выпуклая, что и требовалось доказать.

Очевидно, что градиенты линейных функций, графики которых есть
$\pi_i, i \in J,$ ограничены константой, зависящей только от
множества $D$ и самой функции $ f_1(\cdot,\cdot).$ Верно равенство
$$
    \psi (r(t)) = f_1(r(t)) \;\; \forall t \in [0,T_r].
$$
Ясно, что $\psi (\cdot,\cdot)$ строится однозначно по функции
$f_1(\cdot,\cdot)$ и выбранной кривой $r(\cdot).$ Из сказанного
выше следует, что функция $\psi (\cdot,\cdot)$ есть липшицевая с
константой $L(D,f)$.

Пусть
$$
 \Psi_1(t) = \psi (r(t)) \;\; \forall t \in [0,T_r].
$$
Поскольку
$$
  \vee (\Phi_1 '; 0,T_r)  =  \vee( \Psi_1 ' ; 0,T_r) ,
$$
то из леммы 1 следует, что для некоторых констант $ c_1(D, f_1),
c_2(D, f_1) > 0 $
$$
\vee (\Phi_1 '; 0,T_r) \leq   c_1(D, f_1) + c_2(D, f_1) \vee(r';
0,T_r ).
$$
Если функция $ f_1(\cdot,\cdot)$ не есть дважды непрерывно
дифференцируемая, то ее можно приблизить выпуклой дважды
непрерывно дифференцируемой функцией $ \tilde{f}_1(\cdot,\cdot)$ и
построить соответствующую ей функцию $\tilde{\psi} (\cdot,\cdot)$
так, чтобы значения  функций $\psi (\cdot,\cdot)$, $\tilde{\psi}
(\cdot,\cdot)$ и их производных там, где они существуют, как
угодно мало отличались друг от друга. Но тогда аналогичное будет
верно для функций $\Psi_1(\cdot)$, $\tilde{\Psi}_1(\cdot)$,
построенных по $\psi (\cdot,\cdot), \tilde{\psi} (\cdot,\cdot)$
соответственно, и их производных. Значит написанное выше
неравенство для вариации производных функции $\Psi_1(\cdot)$ верно
для общего случая. Лемма \ref{2lemdifconv} доказана. $\Box$

Из леммы \ref{2lemdifconv} следует, что если $f(\cdot,\cdot)$
представима  в виде разности выпуклых функций, т.е.
$$ f(z) = f_1(z) - f_2(z) \;\;\;\; \forall z \in D,   $$
где $f_i(\cdot,\cdot), i=1,2,$ - выпуклые, то условие
(\ref{difconv3})  c необходимостью выполняется. Действительно, для
произвольной $r(\cdot) \in \wp(D)$ введем обозначения
$$
\Psi_1(t) = f_1(r(t)), \Psi_2(t) = f_2(r(t)) \;\;\; \forall t \in
[0,T_r].
$$
Поскольку \cite{kolmogorovfomin}
$$
\vee (\Phi ';0,T_r) \leq \vee (\Phi_1 '; 0,T_r) + \vee (\Phi_2
';0,T_r)
$$
то, принимая во внимание неравенство (\ref{difconv4}), неравенство
(\ref{difconv3}) с необходимостью выполняется.

Докажем достаточность условия (\ref{difconv3}) для представления
функции $f(\cdot)$ в виде разности выпуклых функций.

Прежде всего покажем, что  для любого $r(\cdot) \in \wp(D)$ и
констант $c_1(D,f), c_2(D,f)>0 $ верно неравенство
$$
\vee (\Phi_n ';0,T_r) \leq c_1(D,f) + c_2(D,f) \vee(r'; 0, T_r),
$$
где $\Phi_n(t) = f_n(r(t))$. Действительно, для любой триангуляции
области $D$ градиенты в точках $r(t_k) \in r(\cdot), t_k \in
[0,T_r],$ линейных функций, графики которых есть грани функции
$f_n(\cdot)$, будут с любой степенью точности $\varepsilon_n$, где
$\varepsilon_n \rightarrow +0$, близки к обобщенным градиентам
функции $f(\cdot)$. Поэтому произвольная конечная сумма
$$
\sum_{i=1}^{N} \mid \Phi_n'(t_i) - \Phi_n'(t_{i+1}) \mid
$$
для больших $n$ будет как угодно мало отличаться от суммы
$$
\sum_{i=1}^{N} \mid \Phi'(t_i) - \Phi'(t_{i+1}) \mid .
$$
А поскольку вариация функции $\Phi_n'(\cdot)$  может только
возрастать при вложенности триангуляций области $D$ при увеличении
$n$, то отсюда и из сказанного выше следует, что
\begin{equation} \vee (\Phi_n';0,T_r) \leq \vee
(\Phi';0,T_r)+\delta(n) \leq c_1(D,f) + c_2(D,f) \vee(r'; 0, T_r),
\label{difconv5}
\end{equation} где $\delta(n) \rightarrow +0 $ при $n \rightarrow
\infty$.

Вариация производных по направлению вдоль произвольного отрезка
суммы выпуклых функций равна сумме вариаций производных этих
выпуклых функций по тому же отрезку. Если будет доказано, что
сумма вариаций  производных всех выпуклых двугранных углов функции
$f_n(\cdot)$ вдоль любого отрезка области $D$  ограничена сверху
константой, независящей от $n$, то отсюда будет следовать, что
ограничена сверху той же константой вариация производной функции
$f_{1,n}(\cdot)$ вдоль произвольного отрезка области $D$. Отсюда
следует равностепенная ограниченность и непрерывность функций
$f_{1,n}(\cdot) - c_{1,n}$. Но тогда по теореме Арцела из
последовательности $f_{1,n}(\cdot) - c_{1,n}$ можно выбрать
подпоследовательность, равномерно сходящуюся на $D$ к выпуклой
функции $f_{1}(\cdot)$. Соответствующая подпоследовательность
последовательности  $f_{2,n}(\cdot) - c_{1,n}$ будет стремиться к
выпуклой функции $f_2(\cdot)$, что означает, что $f(\cdot)$ есть
ПРВ функция.

Пусть условия теоремы выполняются, но $f(\cdot)$ не есть ПРВ
функция. Проделаем следующую процедуру. Путем разбиения множества
$D$ на выпуклые подобласти можно выделить ту подобласть, где
функции $f_{1,n}(\cdot)$ имеют предельное бесконечное значение
вариации производной вдоль некоторых отрезков этой подобласти при
$n \rightarrow \infty$. Действительно, в противном случае из
последовательности функции $f_{1,n}(\cdot) - c_{1,n}$ можно было
бы выбрать сходящуюся подпоследовательность, и $f(\cdot)$ была бы
ПРВ функцией.

Далее разбиваем выделенную подобласть на меньшие области и опять
выделяем ту, где вариация производной функций $f_{1,n}(\cdot)$
вдоль некоторых отрезков неограничена при $n \rightarrow \infty$.
В итоге определяем точку $M$, в произвольной окрестности которой
вариация производной функций $f_{1,n}(\cdot)$ вдоль некоторых
отрезков неограничена при $n \rightarrow \infty$. Без ограничения
общности можно считать, что $M-$ внутренняя точка множества
$\bar{D},$ так как все получаемые в процессе применения алгоритма
функции $-$ равномерно липшицевы и могут быть распространены во
вне множества $\bar{D},$

Берем произвольную окрестность точки $M$ и разбиваем ее на
конечное число секторов. Выбираем произвольный из них, где
вариация производной функций $f_{1,n}(\cdot)$ вдоль некоторых
отрезков неограничена при $n \rightarrow \infty$. Далее выбранный
сектор разбиваем на конечное число секторов и опять выбираем тот
из низ, где вариация производной функций $f_{1,n}(\cdot)$ вдоль
некоторых отрезков неограничена при $n \rightarrow \infty$ и т.д.
Множество выбранных секторов стягивается к некоторому направлению,
определяемому единичным вектором $l$ с вершиной в точке $M$.
Очевидно, что в произвольном секторе $K$ с вершиной с точке $M$,
содержащем вектор $\alpha l$ в $\mbox{int} \, K$, $\alpha
>0$, вариация производной функций $f_{1,n}(\cdot)$ вдоль некоторых
отрезков неограничена при $n \rightarrow \infty$.

Возможны два случая:

a) вариация производных функций $f_{1,n}(\cdot)$ по направлению
$l$ неограничена при $n \rightarrow \infty$;

б) вариация производных функций $f_{1,n}(\cdot)$ по направлению
$\eta,$ перпендикулярном направлению $l$, неограничена при $n
\rightarrow \infty$ .

Сказанное можно перефразировать следующим образом, а именно: сумма
вариаций производных выпуклых двугранных углов функции
$f_{n}(\cdot)$ вдоль указанного направления неограничена при $n
\rightarrow \infty$.

Рассмотрим случай а). Возьмем произвольный сектор $K$, содержащий
вектор $\alpha l$ в $\mbox{int} \, K$, $\alpha >0$. Будем
рассматривать выпуклые двугранные углы функций $f_{1,n}(\cdot)$ из
конуса $K$ для всех $n$.

За счет равномерной липшицевости по $n$ всех двугранных углов
функций $f_{n}(\cdot)$ вариации производных по направлению этих
двугранных углов равномерно непрерывны относительно направления и
$n$.

Для каждого выпуклого $k-$ ого  двугранного  функции
$f_{n}(\cdot)$ выделим отрезок $v_{k,n}$, вариация производной
вдоль которого для $k-$ ого двугранного угла максимальна и равна
$a_{k,n}$. Ясно, что отрезок $v_{k,n}$ должен быть перпендикулярен
проекции на плоскость $XOY$ линии раздела двух граней $k-$ ого
двугранного угла.

Пусть угол наклона отрезков $v_{k,n}$ с направлением $l$ не
превосходит $\pi / 2- \delta$ для некоторого $ \delta>0.$

Путем разбиения сектора $K$ на меньшие секторы, стягивающиеся к
вектору $\alpha l$ и точку $M$, и рассмотрения в каждом из них
своей группы отрезков $v_{k,n}$ для всех значений $k$ и $n$, можно
выделить одну или несколько групп указанных отрезков, каждую из
которых можно пересечь кривой $r(\cdot) \in \wp(D),$ образующей в
точке пересечения с отрезками $v_{k,n}$ угол, не превосходящий
$\pi /2 -\delta, \delta>0.$ Поскольку сектор $K$ произвольный,
содержащий вектор $\alpha l$, то можно рассматривать такие кривые,
для которых $r'(t) \rightarrow -l,$ когда $r(t) \rightarrow M.$
Сама кривая $r(\cdot)$ будет включать в себя отрезки, близкие к
отрезкам $v_{k,n}$.

Если для рассматриваемого случая подгруппа отрезков $\{ v_{k,n}
\}$ существует только одна, то вдоль найденной кривой $ r(\cdot)
\in \wp(D)$ вариация производной суммы выпуклых двугранных углов
  стремится к бесконечности при $n \rightarrow +\infty$.

Кривая $r(\cdot)$, как упоминалось, строится таким образом, чтобы
она включала отрезки, близкие к отрезкам $\{ v_{k,n} \}$. Так как
при выполнении неравенства (\ref{difconv3}) выполняется
неравенство (\ref{difconv5}), а мы нашли кривую $r(\cdot)$, вдоль
которой сумма вариаций производных двугранных углов бесконечна, то
из (\ref{difconv5}) следует, что неограничена вдоль $r(\cdot)$
вариация производной функции $\Phi(\cdot)$. Приходим к
противоречию.

Кроме того, возможен случай, когда у нас есть несколько групп
отрезков $\{ v_{k,n} \}_i,$ для каждой из которых найдется кривая
$r_i(\cdot) \in \wp(D)$, что
$$
\vee (\Phi_n';0,T_{r_i})= c_i, \,\,\,  r_i' (t) \rightarrow_{t
\rightarrow T_{r_i}} -l,
$$
где $T_{r_i}-$  параметр кривой $r_i(\cdot)$ при естественной
параметризации в точке $M$,  а также
$$
 \sum_i \, c_i = \infty.
$$
Тогда кривую $r(\cdot) \in \wp(D),$ вдоль которой сумма вариаций
производных двугранных углов  стремится к бесконечности при $n
\rightarrow +\infty$, будем строить следующим образом.

Кривая $r(\cdot)$ должна содержать достаточное количество $k_i$
отрезков из каждой группы отрезков $\{ v_{k,n} \}_i,$ (либо
близких к ним), чтобы
$$
\vee (\Phi'_n ;t_{r_{i}},t_{r_{i+1}}) =  c_i - \mu_i,
$$
где $t_{r_i}>0-$ значения параметра $t$ для $i$-ой группы отрезков
при естественной параметризации кривой $r_i(\cdot)$, $\mu_i < c_i
-$ малые положительные числа, для которых
$$
\sum_i \, \mu_i < \infty.
$$

Нетрудно видеть, что всегда такую кривую $r(\cdot)$ построить
можно. Она будет состоять из набора  кривых $r_i(\cdot)$.   Для
этого надо осуществить плавный переход от одной кривой
$r_i(\cdot)$ к кривой $r_{i+1}(\cdot),$ не выходя из множества
$\wp(D)$, обеспечив конечную вариацию производной $r'(\cdot)$.
Поскольку $r'_i(t) \rightarrow -l$ при $t \rightarrow T_{r_i}$ для
вех $i$, то подобная процедура осуществима всегда.

Но тогда
$$
\vee (\Phi'_n ; 0,T_r) \geq \sum_i  \vee
(\Phi'_n;t_{r_{i}},t_{r_{i+1}}) =
$$
$$
=\sum_i(c_i-\mu_i)= \sum_i c_i -\sum_i \mu_i  =\infty.
$$

Но тогда, как следует из (\ref{difconv5}), нарушается неравенство
(\ref{difconv3}), которое по предположению достаточности условия
теоремы является верным. Опять приходим к противоречию.

Если вариация производной суммы выпуклых двугранных углов функции
$f_n(\cdot)$ конечна вдоль направления, определяемого вектором
$l$, при любом $n$, то для случая неограниченности при $n
\rightarrow \infty$ вариации производной функции $f_{1,n}(\cdot)$
в произвольно малом секторе с вершиной $M$, содержащем вектор
$\alpha l$, $\alpha >0$, следует, что вариация производной суммы
выпуклых двугранных углов функции $f_n(\cdot)$ бесконечна при $n
\rightarrow \infty$  вдоль направления $\eta $,

Случай б). Все отрезки $v_{k}$ можно разбить на такие группы $\{ m
\}$ отрезков, которые можно пересечь кривой $r_{m}(\cdot) \in
\wp(D),$ для которой
$$
        r'_{m} (\tau) \rightarrow_{\tau \rightarrow T_{r_{m}}} -l,
$$
где $T_{r_{m}}-$ есть параметр кривой $r_{m}(\cdot)$ при
естественной параметризации в точке $M,$  и кривизна кривой $r_{m}
(\cdot)$ стремится к бесконечности при $\tau \rightarrow
T_{r_{m}}.$ Кривая $r_{m}(\cdot)$ пересекает свою группу отрезков
под острыми углами $\alpha_{k_m}$ в точках $\tau_{k_m}$, причем
$\alpha_{k_m} \rightarrow \pi / 2$ при $\tau_{k_m} \rightarrow
T_{r_{m}}$. Ясно, что сказанное всегда выполнимо путем разбиения
множества всех отрезков $v_{k}$ на подмножества с требуемыми
свойствами.

Кроме того, углы $\alpha_{k_m},$ кривые $r_{m}(\cdot)$ и группы
отрезков $\{ v_{k} \}_m$ можно выбрать такими, чтобы предел по $m$
суммы вариаций производных функций $\Phi'_m(\cdot)$ вдоль отрезков
кривых $r_{m}(\cdot)$  был равен бесконечности. В противном случае
функции $f_n(\cdot)$ имели бы ограниченную вариацию вдоль
направления $\eta$ при $n \rightarrow +\infty$ (см. замечание).

Построение кривых $r_{m}(\cdot)$ с неограниченно увеличивающейся
кривизной в точке $M$, для которой
$$
   \lim_{m \rightarrow \infty} \sum_{m}
   \vee (\Phi'_m ; T_{k_m},T_{{k_m+1}})=\infty,
$$
осуществляется аналогичным способом, как и в случае a). Для этого
надо построить кривую $r_{m}(\cdot) \in \wp(D)$ с описанными выше
свойствами, состоящую из достаточного  количества $k_{m}$ отрезков
$\{v_{k}\}$ (либо близких к ним), чтобы
$$
\sum_{k_m} \vee (\Phi'_m ;T_{k_{m}},T_{k_{m}+1}) = c_{m},
$$
и
$$
\lim_m  c_{m} =\infty,
$$
$[ T_{k_{m}}, T_{k_{m}+1} ]$- значение параметра $t$ для отрезка
кривой $r_m (\cdot)$ при ее естественной параметризации. Такие
кривые $r_{m}(\cdot)$ всегда можно построить. При $m \rightarrow
\infty$ кривые $r_{m}$ будут пересекать под острыми углами все
большее число указанных отрезков из произвольно малого сектора,
содержащем вектор $\alpha l$, с вершиной в точке $M$. Кривизны
кривых $r_{m}$ вблизи точки $M$ неограниченно увеличиваются при $m
\rightarrow \infty$.

Но тогда из последовательности кривых $r_m(\cdot)$ можно составить
кривую $r(\cdot) \in \wp(D) $ с конечной суммой вариаций
производной на отрезках $[ T_{k_m},T_{k_m+1} ]$, что для
$\Phi(t)=f(r(t)) $
$$
\sum_{k_m} \vee (\Phi'; T_{k_m},T_{k_m+1}) = \infty.
$$
Отсюда и из (\ref{difconv5}) приходим к противоречию с
(\ref{difconv3}).

Итак, доказано, что при выполнении условия теоремы, сумма вариаций
производных выпуклых двугранных углов функции $f_n(\cdot)$ вдоль
любого отрезка области $D$ при $n \rightarrow \infty$ ограничена
сверху константой, независящей от $n$. Отсюда, как отмечалось
выше, следует, что $f(\cdot)-$ ПРВ функция.

Итак, теорема \ref{difconvexthm1} доказана. $\Box$

{\bf Замечание 1.} {\em Рассуждения с выбором углов $\alpha_{k_m}$
и кривых $r_m(\cdot)$ аналогичны следующим.

Пусть имеем расходящийся ряд
$$
\sum_i \, a_i = \infty, \,\,\,\, a_i >0 \,\,\,\,\, \forall i.
$$
Всегда можно выбрать монотонно убывающую по $i$ последовательность $\{ \beta_i \}, \\
\beta_i \rightarrow_{i \rightarrow \infty} 0,$ чтобы
$$
\lim_{m \rightarrow \infty} \sum_{i=1}^{m} \beta_i \, a_i =
\infty.
$$
Здесь $a_i$ является аналогом вариации производной двугранного
угла вдоль отрезка $v_i$, а $\beta_i$ - аналог косинуса угла,
образуемого кривой $r_i$ с этим отрезком в точке пересечения. }

\vspace{1cm}

\section{Геометрическая интерпретация теоремы 1}

\vspace{0.5cm}

Перефразируем теорему \ref{difconvexthm1}, придав ей более
геометрический характер. Для этого введем понятие поворота кривой
$r(\cdot)$  на графике $\Gamma_f = \{(x,y,z) \in \mathbb{R}^3 \mid
z = f(x , y)\}.$

Рассмотрим на $\Gamma_f$ кривую $ R(t)=(r(t),f(r(t))),$ где
$r(\cdot) \in \wp(D).$  Так как функция $f(\cdot,\cdot)$ есть
липшицевая, то п.в. на $[0,T_r]$ существует производная
$R'(\cdot),$ которую обозначим через $\tau(\cdot)=R'(\cdot).$

{\bf Определение 4.4.1.} {\em  Поворотом кривой $R(\cdot)$ на
многообразии $\Gamma_f$ назовем величину
$$
sup_{ \{t_i\} \subset N_r} \,\, \sum_i \Vert \tau(t_i)/ \Vert \tau
(t_i) \Vert -  \tau(t_{i-1})/ \Vert \tau (t_{i-1}) \Vert \Vert =
O_r.
$$
}

Таким образом, поворот $O_r$ кривой $R(\cdot)$ есть верхняя грань
суммы углов между касательными $\tau(t)$  для $t \in [0,T_r].$
Нетрудно видеть, что для плоской гладкой кривой, параметризованной
естественным образом, величина $ O_r$ равна интегралу
$$
\int^{T_r}_0 \mid k(s) \mid ds,
$$
где $k(s)$ - кривизна рассматриваемой кривой $ r(\cdot)$ в точке
$s \in [0,T_r],$ т.е. совпадает с обычным определением поворота
кривой в точке \cite{pogorelov1} .

{\bf Теорема 2.} {\em Для того, чтобы произвольная липшицевая функция \\
$z \rightarrow f(z) :D \rightarrow \mathbb{R}$ была ПРВ функцией
на выпуклом компактном множестве $ D \in \mathbb{R}^2,$ необходимо
и достаточно, чтобы для любой кривой $r(\cdot) \in \wp(D)$ и любой
ее системы подмножеств существовали константы $c_3(D,f), c_4(D,f)
>0$ такие, что поворот кривой $R(\cdot)$ на $\Gamma_f$ ограничен
сверху неравенством, т.е.
\begin{equation} O_r \leq c_3(D,f) + c_4(D,f) \vee(r';0,T_r)  \;\;
\forall r \in \wp(D). \label{difconv6} \end{equation} }

{\bf Замечание 3}. {\em  Условие теоремы означает следующее. Сумма
поворотов любых участков кривой $R(\cdot)$ на $\Gamma_f$ ограничен
сверху неравенством (\ref{difconv6}). Здесь также, как и выше,
приходится оговаривать правило подсчета поворота кривой
$R(\cdot)$, так как поворот всей кривой $R(t)$, $t \in [0, T_r] $,
может оказаться бесконечной.}

{\bf Доказательство. }{\bf Необходимость}. Пусть $ f(\cdot,\cdot)$
есть ПРВ функция. Покажем, что тогда справедливо неравенство
(\ref{difconv6}). Воспользуемся неравенством, вытекающим из
неравенства треугольника,
$$
\Vert \tau(t_i) / \Vert \tau(t_i) \Vert  - \tau (t_{i-1}) / \Vert
\tau(t_{i-1} \Vert \Vert \leq \Vert r'(t_i) /
\sqrt{ 1+f'^2_t (r(t_i))}  - r'(t_{i-1}) /
\sqrt{ 1+f'^2_t (r(t_{i-1}))} \Vert +
$$
$$
\mid f'_t(r(t_i)) /  \sqrt{ 1+f'^2_t (r(t_i))}- f'_t(r(t_{i-1})) /
\sqrt{ 1+f'^2_t (r(t_{i-1}))} \mid .
$$
Так как $1 \leq  \sqrt{1+f'^2_t (r(t_i))} \leq \sqrt{1+L^2}$ для
всех $ t_i \in [0,T_r]$ , то очевидно, существует такое $c_3 >1,$
для которого верно неравенство

\begin{equation} \Vert r'(t_i) / \sqrt{1+f'^2_t (r(t_i))}-
r'(t_{i-1}) / \sqrt{1+f'^2_t (r(t_{ i-1}))} \Vert \leq c_3 \Vert
r'(t_i) - r'(t_{i-1}) \Vert. \label{difconv7} \end{equation}

Из свойств функции $\theta(x)= x / \sqrt{ 1+x^2}$ следует
неравенство \begin{equation} \mid f'_t (r(t_i)) / \sqrt{ 1+f'^2_t
(r(t_i))} - f'_t (r(t_{i-1})) / \sqrt{1+f'^2_t (r(t_{i-1}))} \mid
\leq \mid f'_t (r(t_i)) - f'_t (r(t_{i-1})) \mid .
\label{difconv8} \end{equation}

Из (\ref{difconv7}) и (\ref{difconv8}) имеем \begin{equation}
\sup_{\{t_i \} \in N_r } \,\, \sum_i \, \Vert \tau(t_i) / \Vert
\tau(t_i) \Vert - \tau (t_{i-1}) / \Vert \tau(t_{i-1}) \Vert \Vert
\leq c_3 (\vee (  r'; 0 ,T_r) + \vee (\Phi' ;  0,T_r) ).
\label{difconv9} \end{equation}

Так как по условию $ f(\cdot,\cdot)-$  ПРВ функция, то согласно
теореме \ref{difconvexthm1}
$$
\vee (\Phi'; 0,T_r) \leq  c_1(D,f) + c_2(D,f) \vee(r';0,T_r),
$$
откуда с учетом (\ref{difconv9}) следует неравенство
(\ref{difconv6}). Необходимость доказана.

{\bf Достаточность}.  Пусть справедливо неравенство
(\ref{difconv6}). Покажем, что $f(\cdot,\cdot)$ - ПРВ функция.
Воспользуемся неравенством
$$
\Vert \tau(t_i) / \Vert \tau( t_i) \Vert  - \tau(t_{i-1}) /  \Vert
\tau (t_{i-1}) \Vert \Vert \geq \mid  f'_t (r(t_i)) / \sqrt{1+f'_t
(r(t_i))} -
$$
\begin{equation} - f'_t (r(t_{i-1})) / \sqrt {1+f'^2_t
(r(t_{i-1}))} \label{difconv10} \end{equation} Из свойств функции
$\theta(x) = x / \sqrt{1+x^2}$ и из $\Vert f'(z) \Vert \leq L$ для
всех $z \in D,$ где производная существует, следует существование
константы $ c_4(L) > 0,$ для которой
$$
\mid f'_t (r(t_i)) / \sqrt{1+f'^2_t (r(t_i))} - f'_t (r(t_{i-1}))
/ \sqrt{1+f'^2_t(r(t_{i-1}))} \geq c_4 \mid f'_t (r(t_i)) - f'_t
(r(t_{i-1})) \mid,
$$
откуда с учетом (\ref{difconv10}) имеем
$$
c_2(D,f) + c_3(D,f) \vee(r';0,T_r) \geq  \sup_{ \{ t_i \} \subset
N_r} \sum_i \Vert \tau(t_i) / \Vert \tau( t_i) \Vert  -
\tau(t_{i-1}) / \Vert \tau (t_{i-1}) \Vert \Vert \geq
$$
$$
\geq c_4 \vee (\Phi';  0,T_r) .
$$
Из теоремы \ref{difconvexthm1} следует, что $f(\cdot)$ - ПРВ
функция. Достаточность доказана.  $\Box$

Возьмем произвольную кривую $r(\cdot) \in \varrho (D)$. Поскольку
вариация производной $\vee (r'; 0,T_r) $ ограничена сверху для
любой кривой $r(\cdot) \in \varrho (D)$, то из Теоремы 1 и Теоремы
2 вытекают следствия.

{\bf Следствие 1}. {\em Для того, чтобы липшицевая функция $z
\rightarrow f(z):D \rightarrow \mathbb{R}$ была представима в виде
разности выпуклых функций, необходимо, чтобы для любой кривой $r
\in \varrho (D)$ и для некоторой константы $c_5(D,f)>0$
выполнялось неравенство
$$
 (\exists c(D,f),) (\forall r \in \varrho(D, f)) \,\,
 \vee (\Phi'; 0,T_r) <  c_5 (D, f),
$$
где $\Phi(t)=f(r(t)) \;\;\;\; \forall t \in [0,T_r]$.
\label{difconvexconl1} }

{\bf Следствие 2}. {\em Для того, чтобы липшицевая функция $z
\rightarrow f(z):D \rightarrow \mathbb{R}$ была представима в виде
разности выпуклых функций, необходимо, чтобы для любой кривой $r
\in \varrho (D)$ и для некоторой константы $c_6(D,f)>0$ поворот
кривой $R(\cdot)$ на $\Gamma_f$ ограничен сверху неравенством,
т.е. $$ O_r \leq c_6(D,f)\;\; \forall r \in \varrho(D). $$ }

В статье \cite{veselyzajicek} приведен пример, подтверждающий
несправедливость достаточности утверждений Следствия 1 и Следствия
2 \cite{proudconvex2}. Фактически авторы статьи
\cite{veselyzajicek} привели пример, показывающий, что кривых
множества $ \varrho(D) $ недостаточно, чтобы проверить
представимость функции в виде разности выпуклых функций. Класс
кривых $\wp(D) $ значительно шире класса $\varrho(D)$.

Прудников Игорь Михайлович

Научный исследовательский центр (НИЦ) Смоленского государственного
медицинского университета, 214018, г. Смоленск, ул. Кирова, 28.


\newpage


\begin{thebibliography}{14}

\bibitem{proudconvex2} Прудников И.М. К вопросу о  представимости функции двух
переменных в виде разности выпуклых функций  Сибирский
математический журнал РАН, 2014. 55, № 6. С.1368-1380.
\bibitem{aleksandrov1} Александров А.Д. О поверхностях, представимых
в виде разности выпуклых функций //  Изв. АН Каз.ССР 1949. N 3. С.
3-20.
\bibitem{aleksandrov2}  Александров А.Д.  Поверхности, представимые
разностями выпуклых функций //  Докл. АН СССР. 1950. Т. 72, №4.
С.613-616.
\bibitem{hiriarturruty} J.B. Hiriart-Urruty, Generalized differentiability,
duality and optimization for problems dealing with differences of
convex functions. In: Convexity and Duality in Optimization
(Groningen, 1984), Lecture Notes in Econom. and Math. Systems,
256, Springer, Berlin-New York, 1985, pp. 37-70.
\bibitem{vesely} L. Vesely, L. Zajicek, Delta -convex mappings between Banach
spaces and applications, Dissertationes Math. (Rozprawy Mat.) 289
(1989). 52pp.
\bibitem{proudconvex1} Прудников И.М. Необходимые и достаточные условия
представимости положительно однородной функции трех переменных в
виде разности выпуклых  функций //  Известия АН РАН  Т.59. N 5,
1992, С.1116-1128.
\bibitem{zalgaller1} Залгаллер В.А. О представимости функции двух переменных
в виде разности выпуклых функций //  Вестник ЛГУ, N 1, 1963,
С.44-45.
\bibitem{Hartman} Hartman P. On functions representable as a difference of convex
functions, Pacific J.Math., 9 (1959), pp.707-713.
\bibitem{demvas} Демьянов В.Ф., Васильев Л.В. Недифференцируемая
оптимизация.  М,: Наука, 1981. 384 С.
\bibitem{clark1} Кларк Ф. Оптимизация и негладкий анализ. Нью-Йорк, 1983.
\bibitem{kolmogorovfomin} Колмогоров А.Н., Фомин С.В. Елементы теории
функций и функционального анализа. М.: Наука, 1976. 544 С..
\bibitem{aleksandrov3} Aleksandrov A.D., Reshetnyak Yu.G. General Theory of
Irregular Curves. Amsterdam: Kluwer Academic Publisher. 1989. 288
P.
\bibitem{pogorelov1} Погорелов А.В. Дифференциальная геометрия. M,: Наука,
1974, 176 С.
\bibitem{veselyzajicek} L.Vesely, L.Zajicek. A non-DC function which is DC
along all convex curves // J. Math. Anal. Appl. 463(2018),
167-175.

\end{thebibliography}
\end{document}